\documentclass[twocolumn,amsmath,amssymb,aps,secnumarabic,nofootinbib,%
    superscriptaddress,floatfix]{revtex4}
\usepackage{amsthm,times,mathptmx,graphicx,pstricks,pst-node}
\newtheorem{theorem}{Theorem}
\newtheorem{lemma}[theorem]{Lemma}
\newtheorem{corollary}[theorem]{Corollary}
\newtheorem{proposition}[theorem]{Proposition}
\newtheorem{question}[theorem]{Question}
\theoremstyle{remark}
\newtheorem*{remark}{Remark}
\newcommand{\Hom}{\operatorname{Hom}}
\newcommand{\conv}{\operatorname{conv}}
\newcommand{\Z}{\mathbb{Z}}
\newcommand{\R}{\mathbb{R}}
\newcommand{\F}{\mathbb{F}}
\newcommand{\HH}{\mathbb{H}}
\renewcommand{\hat}{\widehat}
\renewcommand{\tilde}{\widetilde}
\renewcommand{\d}{\partial}
\newcommand{\st}{\mathrm{st}}
\newcommand{\hg}{{\hat{g}}}
\newcommand{\hS}{{\hat{S}}}
\newcommand{\pol}{{}^\triangle}
\newcommand{\ppol}{{}^\diamondsuit}
\newcommand{\maxfatness}{5.048}
\newcommand{\eqdef}{\stackrel{\mathrm{def}}{=}}

\newcommand{\Eg}{\emph{E.g.}}
\newcommand{\ie}{\emph{i.e.}}
\newcommand{\cf}{\emph{cf.}}
\newcommand{\eatline}{\vspace{-\baselineskip}}

\newenvironment{fullfigure}[2]
    {\begin{figure}[ht]\begin{center}\def\ffa{#1}\def\ffb{#2}}
    {\vspace{\baselineskip}\caption{\ffb.}\label{\ffa}\end{center}\end{figure}}

\newcommand{\fig}[1]{Figure~\ref{#1}}
\newcommand{\thm}[1]{Theorem~\ref{#1}}
\renewcommand{\sec}[1]{Section~\ref{#1}}
\newcommand{\lem}[1]{Lemma~\ref{#1}}
\newcommand{\prop}[1]{Proposition~\ref{#1}}
\newcommand{\eq}[2]{\begin{equation}\label{#1}#2\end{equation}}

\psset{linewidth=.5pt,dash=2.5pt 2.5pt,unit=.25in,arrowsize=4pt 3}
\SpecialCoor

\newgray{gray90}{.90}
\newcommand{\ar}{\pspicture(0,0)(0,0)
    \psline[arrows=->](4.3pt,0)(4.4pt,0)\endpspicture}
\newcommand{\sgpie}{
    \pscustom[linestyle=none,fillstyle=solid,fillcolor=gray90]{
    \psarcn(4.899;30){1.464}{255}{165} \psarcn(4.899;60){1.464}{285}{195}
    \psarcn(4.899;90){1.464}{315}{225} \psarcn(4.899;120){1.464}{345}{255}
    \psarcn(4.899;150){1.464}{15}{285} \psarcn(4.899;210){1.464}{75}{345}
    \psarcn(4.899;240){1.464}{105}{15} \psarcn(4.899;270){1.464}{135}{45}
    \psarcn(4.899;300){1.464}{165}{75} \psarcn(4.899;330){1.464}{195}{105}}
    \psarcn(4.899;30){1.464}{255}{165} \psarcn(4.899;60){1.464}{285}{195}
    \psarcn(4.899;90){1.464}{315}{225} \psarcn(4.899;120){1.464}{345}{255}
    \psarcn(4.899;150){1.464}{15}{285} \psarcn(4.899;210){1.464}{75}{345}
    \psarcn(4.899;240){1.464}{105}{15} \psarcn(4.899;270){1.464}{135}{45}
    \psarcn(4.899;300){1.464}{165}{75} \psarcn(4.899;330){1.464}{195}{105}
    \qdisk(4;165){.1} \qdisk(4;195){.1} \qdisk(4;225){.1} \qdisk(4;255){.1}
    \qdisk(4;285){.1} \qdisk(4;315){.1} \qdisk(4;345){.1} \qdisk(4;15){.1}
    \qdisk(4;45){.1} \qdisk(4;75){.1} \qdisk(4;105){.1} \qdisk(4;135){.1}
    \rput{120}(3.435;30){\ar} \rput{150}(3.435;60){\ar}
    \rput{180}(3.435;90){\ar} \rput{210}(3.435;120){\ar}
    \rput{240}(3.435;150){\ar} \rput{300}(3.435;210){\ar}
    \rput{330}(3.435;240){\ar} \rput{0}(3.435;270){\ar}
    \rput{30}(3.435;300){\ar} \rput{60}(3.435;330){\ar}
    \uput[180](4;180){\large $\vdots$}
    \uput[0](4;0){\large $\vdots$}}

\begin{document}
\title{Fat 4-polytopes and fatter 3-spheres}

\author{David Eppstein}
\email[Email: ]{eppstein@ics.uci.edu}
\thanks{Supported in part by NSF grant CCR \#9912338}
\affiliation{Department of Information and Computer Science,
    University of California, Irvine, CA 92697}
\author{Greg Kuperberg}
\email[Email: ]{greg@math.ucdavis.edu}
\thanks{Supported by NSF grant DMS \#0072342}
\affiliation{Department of Mathematics,
    University of California, Davis, CA 95616}
\author{G\"unter M. Ziegler}
\email[Email: ]{ziegler@math.tu-berlin.de}
\thanks{Supported by a DFG Leibniz grant}
\affiliation{Institut f\"ur Mathematik, MA 6-2,
    Technische Universit\"at Berlin, D-10623 Berlin, Germany}

\begin{abstract}
We introduce the \emph{fatness} parameter of a $4$-dimensional polytope $P$,
defined as $\phi(P)=(f_1+f_2)/(f_0+f_3)$.  It arises in an important open
problem in $4$-dimensional combinatorial geometry: Is the fatness of convex
$4$-polytopes bounded?

We describe and analyze a hyperbolic geometry construction that produces
$4$-polytopes with fatness $\phi(P)>\maxfatness$, as well as the first infinite
family of $2$-simple, $2$-simplicial $4$-polytopes. Moreover, using a
construction via finite covering spaces of surfaces, we show that fatness is
not bounded for the more general class of strongly regular CW decompositions of
the $3$-sphere.
\end{abstract}

\maketitle

\section{Introduction}
\label{s:intro}

The characterization of the set ${\cal F}_3$ of $f$-vectors of convex 3-dimensional
polytopes (from 1906, due to Steinitz \cite{Steinitz:eulerschen}) is well-known
and explicit, with a simple proof: An integer vector $(f_0,f_1,f_2)$ is the
$f$-vector of a $3$-polytope if and only if it satisfies
\begin{itemize}\setlength{\itemsep}{0pt}
\item $f_1=f_0+f_2-2$ (the Euler equation),
\item $f_2\le 2f_0-4$ (with equality for simplicial polytopes), and
\item $f_0\le 2f_2-4$ (with equality for simple polytopes).
\end{itemize}
(Recall that by the definition of the $f$-vector, $f_k$ is the number
of $k$-faces of the polytope.)
This simple result is interesting for several reasons:
\begin{itemize}\setlength{\itemsep}{0pt}
\item The set of $f$-vectors is the set of \emph{all} the integer points in a
closed $2$-dimensional \emph{polyhedral cone} (whose apex is the $f$-vector
$f(\Delta_3)=(4,6,4)$ of a $3$-dimensional simplex). In particular, it is
\emph{convex} in the sense that ${\cal F}_3=\conv({\cal F}_3)\cap \Z^3$.
\item The same characterization holds for convex $3$-polytopes (geometric
objects), more generally for strongly regular CW $2$-spheres (topological
objects), and yet more generally for Eulerian lattices of length $4$
(combinatorial objects \cite{Stanley:eulerian}).
\end{itemize}

In contrast to this explicit and complete description of ${\cal F}_3$, our
knowledge of the set ${\cal F}_4$ of $f$-vectors of (convex) $4$-polytopes (see
Bayer \cite{Bayer:extended} and H\"oppner and Ziegler \cite{HZ:census}) is very
incomplete. We know that the set ${\cal F}_4$ of all $f$-vectors of
4-dimensional polytopes has no similarly simple description. In particular, the
convex hull of ${\cal F}_4$ is not a cone, it is not a closed set, and not all
integer points in the convex hull are $f$-vectors. Also, the $3$-dimensional
cone with apex $f(\Delta_4)$ spanned by ${\cal F}_4$ is not closed, and its
closure may not be polyhedral.

Only the two extreme cases of simplicial and of simple $4$-polytopes (or
$3$-spheres) are well-understood.   Their $f$-vectors correspond to faces of
the convex hull of ${\cal F}_4$, defined by the valid inequalities $f_2\ge
2f_3$ and $f_1\ge 2f_0$, and the $g$-Theorem, proved for $4$-polytopes by
Barnette \cite{Barnette:inequalities} and for $3$-spheres by Walkup
\cite{Walkup:lower}, provides \emph{complete} characterizations of their
$f$-vectors. (The $g$-Theorem for general simplicial polytopes was famously
conjectured by McMullen \cite{McMullen:simplicial} and proved by Billera and
Lee \cite{BL:sufficiency} and Stanley \cite{Stanley:simplicial}.  See
\cite[\S8.6]{Ziegler:gtm} for a review.)

But we have no similarly complete picture of other extremal types of
$4$-polytopes. In particular, we cannot currently answer the following key
question: Is there a constant $c$ such that all  $4$-dimensional convex
polytopes $P$ satisfy the inequality
\[
f_1(P)+f_2(P)\;\le\;c\,(f_0(P)+f_3(P))?
\]
To study this question, we introduce the \emph{fatness} parameter
\[
\phi(P)\;\eqdef\;\frac{f_1(P)+f_2(P)}{f_0(P)+f_3(P)}
\]
of a $4$-polytope $P$.  We would like to know whether fatness is bounded.

For example, the $4$-simplex has fatness $2$, while the $4$-cube and the
$4$-cross polytope have fatness $\frac{56}{24}=\frac73$. More generally, if $P$
is simple, then we can substitute the Dehn-Sommerville relations
\[
f_2(P)=f_1(P)+f_3(P)-f_0(P)\qquad f_1(P)=2f_0(P)
\]
into the formula for fatness, yielding
\[
\phi(P)=\frac{f_1(P)+f_2(P)}{f_0(P)+f_3(P)}
    = \frac{3f_0(P)+f_3(P)}{f_0(P)+f_3(P)}<3.
\]

Since every $4$-polytope and its dual have the same fatness, the same upper
bound holds for simplicial $4$-polytopes. On the other hand, the ``neighborly
cubical'' $4$-polytopes of Joswig and Ziegler \cite{JZ:cubical} have
$f$-vectors
\[
(4,2n,3n-6,n-2)\cdot 2^{n-2},
\]
and thus fatness
\[
\phi=\frac{5n-6}{n+2} \to 5.
\]
In particular, the construction of these polytopes disproved the conjectured
flag-vector inequalities of Bayer \cite[pp.~145, 149]{Bayer:extended} and
Billera and Ehrenborg \cite[p.~109]{HZ:census}.

The main results of this paper are two lower bounds on fatness:

\begin{theorem} There are convex $4$-polytopes $P$ with fatness
$\phi(P)>\maxfatness$. \label{th:fat} \end{theorem}

\begin{theorem} The fatness of cellulated $3$-spheres is not bounded.
A $3$-sphere $S$ with $N$ vertices may have fatness as high as
$\phi(S) = \Omega(N^{1/12})$.
\label{th:fatter} \end{theorem}

We will prove \thm{th:fat} in \sec{s:poly} and \thm{th:fatter} in
\sec{s:spheres}, and present a number of related results along the
way.

\section{Conventions}
\label{s:conven}

Let $X$ be a finite
CW complex. If $X$ is identified with a manifold $M$, it is also called a
\emph{cellulation} of $M$. The complex $X$ is \emph{regular} if all closed
cells are embedded \cite[\S38]{Munkres:algtop}.
If $X$ is regular, we define it to
be \emph{strongly regular} if in addition the intersection of any
two closed cells
is a cell.  For example, every simplicial complex is a strongly regular CW
complex.  The complex $X$ is \emph{perfect} if the boundary maps of its chain
complex vanish.  (A non-zero-dimensional perfect complex is never regular.)

The \emph{$f$-vector} of a cellulation $X$, denoted $f(X)=(f_0,f_1,\ldots)$,
counts the number of cells in each dimension:  $f_0(X)$ is the number of
vertices, $f_1(X)$ is the number of edges, etc.

If $X$ is $2$-dimensional, we define its fatness as
\[
\phi(X) \eqdef \frac{f_1(X)}{f_0(X) + f_2(X)}.
\]
If $X$ is $3$-dimensional, we define its fatness $X$ as
\[
\phi(X) \eqdef \frac{f_1(X) + f_2(X)}{f_0(X) + f_3(X)}.
\]
If $P$ is a convex $d$-polytope, then its $f$-vector $f(P)$ is defined to be
the $f$-vector of its boundary complex, which is a strongly regular
$(d-1)$-sphere. If $d=4$ we can thus consider $\phi(P)$, the fatness of $P$.
The faces of $P$ of dimension 0 and 1 are called \emph{vertices} and
\emph{edges} while the faces of dimension $d-1$ and $d-2$ are called
\emph{facets} and \emph{ridges}. We extend this terminology to general
cellulations of $(d-1)$-manifolds.

The \emph{flag vector} of a regular cellulation $X$, and likewise the flag
vector of a polytope $P$, counts the number of nested sequences of cells with
prespecified dimensions.  For example, $f_{013}(X)$ is the number triples
consisting of a $3$-cell of $X$, an edge of the $3$-cell, and a vertex of the
edge; if $P$ is a 4-cube, then
$$f_{013}(P) = 192.$$

A convex polytope $P$ is \emph{simplicial} if each facet of $P$ is a simplex.
It is \emph{simple} if its polar dual $P\pol$
is simplicial, or equivalently if the cone of each vertex matches
that of a simplex.

\section{4-Polytopes}
\label{s:poly}

In this section we construct families of $4$-polytopes with several
interesting properties:

\begin{itemize}\setlength{\itemsep}{0pt}
\item They are the first known infinite families of {\em $2$-simple,
$2$-simplicial $4$-polytopes}, that is, polytopes in which all $2$-faces and
all dual $2$-faces are triangles (all edges are ``co-simple''). \Eg, Bayer
\cite{Bayer:extended} says that it would be interesting to have an infinite
family.

As far as we know, there were only six such polytopes previously known: the
simplex, the hypersimplex (the set of points in $[0,1]^4$ with coordinate sum
between $1$ and $2$), the dual of the hypersimplex, the $24$-cell, and a gluing
of two hypersimplices (due to Braden \cite{Braden:glued}), and the dual of the
gluing. (There are claims in Gr\"unbaum
\cite[p.~82, resp. p.~170]{Grunbaum:convex} that results of Shephard, resp.\
Perles and Shephard, imply the existence of infinitely many $2$-simple
$2$-simplicial $4$-polytopes.  Both claims appear to be incorrect.)

\item They are the fattest known convex $4$-polytopes.

\item They yield finite packings of (not necessarily congruent) spheres in
$\R^3$ with slightly higher average kissing numbers than previously known
examples \cite{Kuperberg:kissing}.
\end{itemize}

Let $Q\subset\R^4$ be a $4$-polytope that contains the origin in its interior.
If an edge $e$ of $Q$ is tangent to the unit sphere $S^3\subset\R^4$ at a point
$t$, then the corresponding ridge ($2$-dimensional face) $F=e\ppol$ of the
polar dual $P=Q\pol$ is also tangent to $S^3$ at $t$.  (Recall that the
polar dual is defined as
$$Q\pol \eqdef \{p | p \cdot q \le 1 \forall q \in Q\}.)$$
Furthermore, the affine
hulls of $e$ and $F$ form orthogonal complements in the tangent space of $S^3$,
so the convex hull $\conv(e\cup F)$ is an orthogonal bipyramid tangent to $S^3$
(\cf, Schulte \cite[Thm.~1]{Schulte:steinitz}).

We will construct polytopes $E$ by what we call the \emph{$E$-construction}.
This means that they are convex hulls
$$E\eqdef\conv(Q\cup P),$$
where $Q$ is a simplicial $4$-dimensional polytope whose edges are tangent to
the unit $3$-sphere $S^3$, and $P$ is the polar dual of $Q$.  Thus $P$ is
simple and its ridges are tangent to $S^3$. (We call $Q$ \emph{edge-tangent}
and $P$ \emph{ridge-tangent}.)

\begin{proposition}
If $P$ is a simple, ridge-tangent $4$-polytope, then the $4$-polytope
$E=\conv(P\cup Q)$ produced by the $E$-construction is $2$-simple and
$2$-simplicial, with $f$-vector
\[
 f(E) = (f_2(P), 6f_0(P), 6f_0(P), f_2(P)),
\]
and fatness
\[
 \phi(E) = 6\,\frac{f_0(P)}{f_2(P)}.
\]
\label{p:e-construction} \end{proposition}

\begin{proof}
Another way to view the $E$-construction is that $E$ is produced from $P$ by
adding the vertices of $Q$ sequentially. At each step, we cap a facet of $P$
with a pyramid whose apex is a vertex of $Q$. Thus the new facets consist of
pyramids over the ridges of $P$, where two pyramids with the same base
(appearing in different steps) lie in the same hyperplane (tangent to $S^3$),
and together form a bipyramid.  The facets of the final polytope $E$ are
orthogonal bipyramids over the ridges of $P$ and are tangent to $S^3$.  Since
the 2-faces of $E$ are pyramids over the edges of $P$, $E$ is $2$-simplicial.

The polytope $E$ is $2$-simple if and only if each edge is \emph{co-simple},
\ie, contained in exactly three facets of $E$.  The iterative construction of
$E$ shows that it has two types of edges: (i) edges of $P$, which are co-simple
in $E$ if and only if they are co-simple in $P$, and (ii) edges formed by
adding pyramids, which are co-simple if and only if the facets of $P$ are
simple.  Since $P$ is simple, its facets are simple and its edges are
co-simple, so $E$ is then $2$-simple.

This combinatorial description of $E$ yields an expression for the $f$-vector
of $E$ in terms of the flag vector
\cite{BJ:dehn,Bayer:extended,HZ:census} of $P$.  Since the facets
of $E$ are bipyramids over the ridges of $P$, the following identities hold:
\begin{align*}
f_3(E) &= f_2(P) & f_2(E) &= f_{13}(P) \\
f_1(E) &= f_1(P)+f_{03}(P) & f_0(E) &= f_0(P)+f_3(P).
\end{align*}
Since $P$ is simple,
\[
f_{03}(P)=4f_0(P) \qquad f_{13}(P)=3f_1(P).
\]
These identities together with Euler's equation and $f_1(P)=2f_0(P)$ imply the
proposition.
\end{proof}

The $f$-vector of $P$ also satisfies
\[
f_0(P)-f_1(P)+f_2(P)-f_3(P)=0 \qquad f_1(P)=2f_0(P),
\]
in the second case because $P$ is simple, so the fatness of $E$ can also be
written
\[
\phi(E) = 6\,\Big(1-\frac{f_3(P)}{f_2(P)}\Big) =
6\,\Big(1-\frac{f_0(Q)}{f_1(Q)}\Big).
\]
Thus maximizing the fatness of $E$ is equivalent to maximizing the ridge-facet
ratio $f_2(P)/f_3(P)$, or the average degree of the graph of $Q$.  It also
shows that the $E$-construction cannot achieve a fatness of $6$ or more.

In light of \prop{p:e-construction}, we would like to construct edge-tangent
simplicial $4$-polytopes. Regular simplicial $4$-polytopes (suitably scaled)
provide three obvious examples: the $4$-simplex $\Delta_4$, the cross polytope
$C_4\pol$, and the $600$-cell.  From these, the $E$-construction produces the
dual of the hypersimplex, the 24-cell, and a new $2$-simple, $2$-simplicial
polytope with $f$-vector $(720,3600,3600,720)$ and fatness $5$,
whose facets are bipyramids over pentagons.

We will construct new edge-tangent simplicial $4$-polytopes by gluing together
(not necessarily simplicial) edge-tangent $4$-polytopes, called \emph{atoms},
to form \emph{compounds}.  We must position the polytopes so that their
facets match, they remain edge-tangent, and the resulting compound is convex.
It will be very useful to interpret the interior of the $4$-dimensional unit
ball as the Klein model of hyperbolic $4$-space $\HH^4$, with $S^3$ the sphere
at infinity.  (See Iversen \cite{Iverson:hyperbolic} and Thurston
\cite[Chap.~2]{Thurston:geometry} for introductions to hyperbolic geometry.) In
particular, Euclidean lines are straight in the Klein model, Euclidean
subspaces are flat, and hence any intersection of a convex polytope with
$\HH^4$ is a convex (hyperbolic) polyhedron.  Even though the Klein model
respects convexity, it does not respect angles. However, angles and convexity
are preserved under hyperbolic isometries.  There are enough isometries to
favorably position certain $4$-polytopes to produce convex compounds.

If a polytope $Q$ is edge-tangent to $S^3$, then it is hyperbolically
hyperideal:  Not only its vertices, but also its edges, lie beyond the sphere
at infinity, except for the tangency point of each edge.  Nonetheless portions
of its facets and ridges lie in the finite hyperbolic realm.  As a hyperbolic
object the polytope $Q$ (more precisely, $Q\cap \HH^4$) is convex and has flat
facets. The ridge $r$ between any two adjacent facets has a well-defined
hyperbolic dihedral angle, which is strictly between $0$ and $\pi$ if (as in
our situation) the ridge properly intersects $\HH^4$.  To compute this angle we
can intersect $r$ at any point $t$ with any hyperplane $R$ that contains the
(hyperbolic) orthogonal complement to $r$ at $t$.  We let $t$ be the tangent
point of any edge $e$ of the ridge, and let $R$ be the hyperplane perpendicular
to $e$.

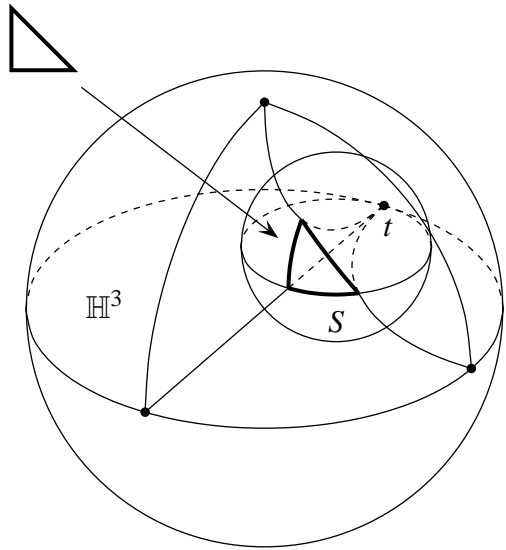
\begin{fullfigure}{f:hyplink}{A cone emanating from an ideal point $t$ in
    $\HH^3$ in the Poincar\'e model, and a horosphere $S$ incident to $t$.  The
    link of $t$ (here a right isosceles triangle) inherits Euclidean geometry
    from $S$}
\begin{large} \pspicture(-5.5,-5.2)(5.2,6.5)
\pscircle(0,0){5}
\pscurve[linestyle=dashed](2.5,2.165)(1.294,2.415)(0,2.5)(-1.294,2.415)
    (-2.5,2.165)(-3.536,1.768)(-4.33,1.25)(-4.83,0.647)(-5,0)
\pscurve[linestyle=dashed](2.5,2.165)(3.536,1.768)(4.33,1.25)(4.83,0.647)(5,0)
\pscurve(-5,0)(-4.83,-0.647)(-4.33,-1.25)(-3.536,-1.768)(-2.5,-2.165)
    (-1.294,-2.415)(0,-2.5)(1.294,-2.415)(2.5,-2.165)(3.536,-1.768)(4.33,-1.25)
    (4.83,-0.647)(5,0)
\pscircle(1.5,1.299){2}
\pscurve[linestyle=dashed](2.5,2.165)(2.018,2.265)(1.5,2.299)(0.982,2.265)
    (0.5,2.165)(0.086,2.006)(-0.232,1.799)(-0.432,1.558)(-0.5,1.299)
\pscurve[linestyle=dashed](2.5,2.165)(3.032,1.942)(3.379,1.641)(3.5,1.299)
\pscurve(-0.5,1.299)(-0.379,0.957)(-0.032,0.656)(0.5,0.433)(1.153,0.314)
    (1.847,0.314)(2.5,0.433)(3.032,0.656)(3.379,0.957)(3.5,1.299)
\pscurve[linestyle=dashed](2.5,2.165)(2.121,1.887)(1.751,1.715)(1.398,1.654)
    (1.07,1.705)(0.776,1.866)
\pscurve(0,4.33)(0.033,3.66)(0.13,3.064)(0.289,2.559)(0.507,2.156)(0.776,1.866)
\pscurve[linestyle=dashed](2.5,2.165)(2.171,1.822)(1.95,1.459)(1.841,1.082)
    (1.848,0.702)(1.97,0.327)
\pscurve(4.33,-1.25)(3.665,-1.02)(3.082,-0.74)(2.597,-0.416)(2.224,-0.057)
    (1.97,0.327)
\psline[linestyle=dashed](2.5,2.165)(0.5,0.433)
\psline(0.5,0.433)(-2.5,-2.165)
\pscurve(-2.5,-2.165)(-2.415,-0.971)(-2.165,0.29)(-1.768,1.531)(-1.25,2.667)
    (-0.647,3.622)(-0,4.33)
\pscurve(4.33,-1.25)(4.183,-0.087)(3.75,1.083)(3.062,2.178)(2.165,3.125)
    (1.121,3.859)(0,4.33)
\pscurve[linewidth=1.5pt](0.5,0.433)(0.512,0.706)(0.546,0.992)(0.602,1.285)
    (0.68,1.579)(0.776,1.866)
\pscurve[linewidth=1.5pt](0.5,0.433)(0.774,0.367)(1.065,0.323)(1.366,0.301)
    (1.67,0.303)(1.97,0.327)
\pscurve[linewidth=1.5pt](0.776,1.866)(0.973,1.581)(1.22,1.244)(1.491,0.894)
    (1.751,0.578)(1.97,0.327)
\qdisk(4.33,-1.25){1.8pt}
\qdisk(-2.5,-2.165){1.8pt}
\qdisk(2.5,2.165){1.8pt}
\qdisk(0,4.33){1.8pt}
\pspolygon[linewidth=1.5pt](-5.3,5)(-5.3,6.3)(-4,5)
\pcline[nodesepA=1.3,nodesepB=1,arrows=->](-4.867,5.433)(1.082,0.875)
\uput[-75](2.5,2.165){$t$}
\rput(1.5,-0.26){$S$}
\rput(-3.348,0.1){$\HH^3$}
\endpspicture
\end{large} \eatline \end{fullfigure}

Within the hyperbolic geometry of $R \cong \HH^3$, every line emanating from
the ideal point $t$ is orthogonal to any horosphere incident to $t$. Thus the
link of the edge $e$ of $Q$ is the intersection of $Q\cap R$ with a
sufficiently small horosphere $S$ at $t$.  Since horospheres have flat
Euclidean geometry \cite[p.~61]{Thurston:geometry}, the link $Q \cap S$ is a
Euclidean polygon.  Its edges correspond to the facets of $Q$ that contain $e$,
and its vertices correspond to the ridges of $Q$ that contain $e$. Thus the
dihedral angle of a ridge $r$ of $Q$ equals to the Euclidean angle of the
vertex $r \cap S$ of the Euclidean polygon $Q \cap S$. This is easier to see in
the Poincar\'e model of hyperbolic space, because it respects angles, than in
the Klein model. \fig{f:hyplink} shows an example.

To summarize:

\begin{lemma}
A compound of two or more polytopes is convex if and only if each ridge has
hyperbolic dihedral angle less than $\pi$, or equivalently, iff each edge
link is a convex Euclidean polygon.
\label{l:convex} \end{lemma}

Compounds can also have interior ridges with total dihedral
angle exactly $2\pi$.  But since all atoms of a compound are
edge-tangent, compounds do not have any interior edges or vertices.

If $Q$ is a regular polytope, then $Q \cap S$ is a regular polygon.
The following lemma is then immediate:

\begin{lemma}
If $Q$ is a regular, edge-tangent, simplicial $4$-polytope,
then in the hyperbolic metric of the Klein model, its dihedral angles
are $\pi/3$ (for the simplex), $\pi/2$ (for the cross polytope), and
$3\pi/5$ (for the $600$-cell).
\label{l:regang} \end{lemma}

A hyperideal hyperbolic object, even if it is an edge-tangent convex
polytope, can be unfavorably positioned so that it is unbounded as a Euclidean
object (\cf, Schulte \cite[p.~508]{Schulte:steinitz}).  Fortunately
there is always a bounded position as well:

\begin{lemma} Let be $Q$ an edge-tangent, convex polytope in $\R^d$ whose
points of tangency with $S^{d-1}$ do not lie in a hyperplane. Then there is a
hyperbolic isometry $h$ (extended to all of $\R^d$) such that $h(Q)$ is
bounded.
\label{l:bounded} \end{lemma}
\begin{proof} Let $p$ lie in the interior of the convex hull of the edge
tangencies and let $f$ be any hyperbolic motion that moves $p$ to the Euclidean
origin in $\R^d$.  Since the convex hull $K$ of the edge tangencies of $h(Q)$
contains the origin, $K\pol$ is a bounded polytope that circumscribes
$S^{d-1}$. Since $K\pol$ is facet-tangent where $h(Q)$ is edge-tangent, $h(Q)
\subset K\pol$.
\end{proof}

In the following we discuss three classes of edge-tangent simplicial convex
$4$-polytopes that are obtained by gluing in the Klein model: Compounds of
simplices, then simplices and cross polytopes, and finally compounds from cut
$600$-cells.  There are yet other edge-tangent compounds involving cross
polytopes cut in half (\ie, pyramids over octahedra), $24$-cells, and
hypersimplices as atoms, but we will not discuss these here.

\subsection{Compounds of simplices}

In this section we classify compounds whose atoms are simplices. This includes
all \emph{stacked} polytopes, which are simplicial polytopes that decompose as
a union of simplices without any interior faces other than facets. However
compounds of simplices are a larger class, since they may have interior ridges.

\begin{lemma}
Any edge-tangent $d$-simplex is hyperbolically regular.
\label{l:hypreg} \end{lemma}

\begin{proof} The proof is by induction on $d$, starting from the case $d=2$, where the
three tangency points define an ideal triangle in $\HH^2$. All ideal triangles
are congruent \cite[p.~83]{Thurston:geometry}.  Since edge-tangent triangles
are the polar duals of ideal triangles, they are all equivalent as well. 

If $d>2$, let $B$ be a general edge-tangent $d$-simplex. On the one hand, there
exists an edge-tangent simplex $A$ which is regular both in hyperbolic geometry
and Euclidean geometry.   On the other hand, given the position of $d$ of the 
vertices, there are at most two choices for the last vertex that  produce an
edge-tangent simplex, one on each side of the hyperplane spanned by the first
$d$. By induction there exists an isometry that takes a face of $B$ to a face
of $A$ and the remaining vertex to the same side.  The edge-tangent constraint
implies that this isometry takes the last vertex of $B$ to the last vertex of
$A$ as well.
\end{proof}

\begin{proposition}\label{p:simplices}
There are only three possible edge-tangent compounds of $4$-simplices:
\begin{itemize}\setlength{\itemsep}{0pt}
\item the regular simplex,
\item the bipyramid (a compound of two simplices that share a facet), and
\item the join of a triangle and a hexagon (a compound of six
  simplices that share a ridge).
\end{itemize}
\end{proposition}

\begin{proof} \fig{f:tjewels} shows all strictly convex polygons with
unit-length edges tiled by unit equilateral triangles, or \emph{triangle
jewels}.  Since the atoms of an edge-tangent compound of simplices are
edge-tangent, they are hyperbolically regular by \lem{l:hypreg}, and their edge
links are equilateral triangles. Thus every edge link of a compound of
simplices is a triangle jewel.  Any three $4$-simplices in a chain in such a
compound share a ridge. In order to create an edge link matching
\fig{f:tjewels}, they must extend to a ring of six simplices around the same
ridge. Adding any further simplex to these six would create an edge link in the
form of a triangle surrounded by three other triangles, which does not appear
in \fig{f:tjewels}.
\end{proof}

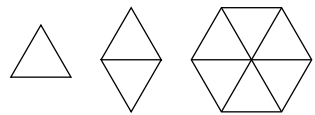
\begin{fullfigure}{f:tjewels}{The 3 possible edge links of edge-tangent
    compounds of $4$-simplices}
\psset{unit=.8cm} \pspicture(-2.5,-1)(3.5,1)
\rput(-1.5,0){\pspicture(-.5,-.6)(.5,.6) 
\pspolygon(0,.577)(-.5,-.289)(.5,-.289)
\endpspicture}
\rput(0,0){\pspicture(-.5,-1)(.5,1) 
\psline(-.5,0)(0,.866)(.5,0)(-.5,0)(0,-.866)(.5,0)
\endpspicture}
\rput(2,0){\pspicture(-1,-1)(1,1) 
\pspolygon(1;0)(1;60)(1;120)(1;180)(1;240)(1;300)
\psline(1;0)(1;180) \psline(1;60)(1;240) \psline(1;120)(1;300)
\endpspicture}
\endpspicture \eatline \end{fullfigure}

Two of the $E$-polytopes produced by \prop{p:simplices} were previously known.
If $Q$ is the simplex, then $E$ is dual to the hypersimplex.  If $Q$ is the
bipyramid, then $E$ is dual to Braden's glued hypersimplex.  However, if $Q$ is
the six-simplex compound (dual to the product of a triangle and a hexagon),
then $E$ is a new $2$-simple, $2$-simplicial polytope with $f$-vector
$(27,108,108,27)$.

\prop{p:simplices} also implies an interesting impossibility
result.

\begin{corollary}
No stacked $4$-polytope with more than $6$ vertices is edge-tangent.
\end{corollary}

See Schulte~\cite[Thm.~3]{Schulte:steinitz} for the first examples of
polytopes that have no edge-tangent realization.

\subsection{Compounds of simplices and cross polytopes}

Next, we consider compounds of simplices and regular cross polytopes. The edge
link of any convex compound of these two types of polytopes must be one of the
eleven strictly convex polygons tiled by unit triangles and squares, or
\emph{square-triangle jewels} (\fig{f:stjewels}).  See Malkevitch
\cite{Malkevitch:tiling} and Waite \cite{Waite:combining} for work on convex
compounds of these shapes relaxing the requirement of strict convexity.

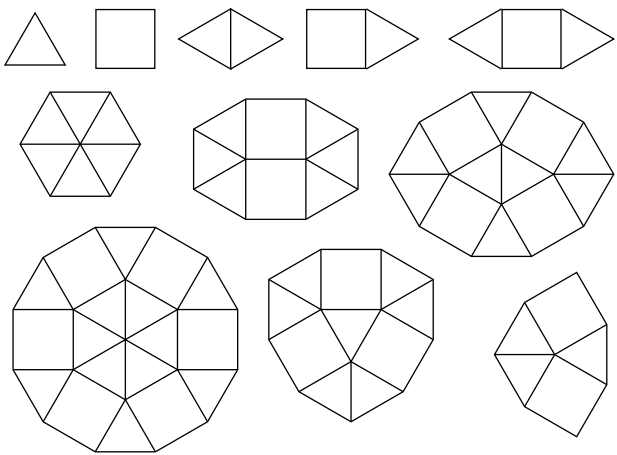
\begin{fullfigure}{f:stjewels}{The 11 possible edge links of edge-tangent
    compounds of $4$-simplices and cross polytopes}
\psset{unit=.8cm} \pspicture(-.75,.75)(10,8.75)
\rput(0,8){\pspicture(-.5,-.5)(.5,.5) 
\pspolygon(0,.433)(-.5,-.433)(.5,-.433) \endpspicture}
\rput(1.5,8){\pspicture(-.5,-.5)(.5,.5) 
\psframe(-.5,-.5)(.5,.5) \endpspicture}
\rput(3.25,8){\pspicture(-1,-.5)(1,.5) 
\psline(0,-.5)(.866,0)(0,.5)(0,-.5)(-.866,0)(0,.5) \endpspicture}
\rput(5.5,8){\pspicture(-.5,-.5)(1.5,.5) 
\psframe(-.5,-.5)(.5,.5) \psline(.5,.5)(1.366,0)(.5,-.5)
\endpspicture}
\rput(8.25,8){\pspicture(-1.5,-.5)(1.5,.5) 
\psframe(-.5,-.5)(.5,.5) \psline(.5,.5)(1.366,0)(.5,-.5)
\psline(-.5,.5)(-1.366,0)(-.5,-.5) \endpspicture}
\rput(.75,6.25){\pspicture(-1,-1)(1,1) 
\pspolygon(1;0)(1;60)(1;120)(1;180)(1;240)(1;300)
\psline(1;0)(1;180) \psline(1;60)(1;240) \psline(1;120)(1;300)
\endpspicture}
\rput(4,6){\pspicture(-1.5,-1)(1.5,1) 
\pspolygon(.5,1)(-.5,1)(-1.366,.5)(-1.366,-.5)(-.5,-1)
    (.5,-1)(1.366,-.5)(1.366,.5)
\psline(-.5,-1)(-.5,1) \psline(.5,-1)(.5,1)
\psline(1.366,.5)(.5,0)(-.5,0)(-1.366,.5)
\psline(1.366,-.5)(.5,0) \psline(-.5,0)(-1.366,-.5)
\endpspicture}
\rput(7.75,5.75){\pspicture(-1.5,-1)(1.5,1) 
\pspolygon(.5,1.366)(-.5,1.366)(-1.366,.866)(-1.866,0)(-1.366,-.866)
    (-.5,-1.366)(.5,-1.366)(1.366,-.866)(1.866,0)(1.366,.866)
\pspolygon(0,-.5)(.866,0)(0,.5)(-.866,0)
\psline(1.866,0)(.866,0) \psline(-1.866,0)(-.866,0)
\psline(1.366,.866)(.866,0)(1.366,-.866)
\psline(-1.366,.866)(-.866,0)(-1.366,-.866)
\psline(.5,1.366)(0,.5)(-.5,1.366) \psline(.5,-1.366)(0,-.5)(-.5,-1.366)
\psline(0,.5)(0,-.5)
\endpspicture}
\rput(1.5,3){\pspicture(-2,-2)(2,2) 
\pspolygon(1.932;15)(1.932;45)(1.932;75)(1.932;105)(1.932;135)(1.932;165)
    (1.932;195)(1.932;225)(1.932;255)(1.932;285)(1.932;315)(1.932;345)
\pspolygon(1;30)(1;90)(1;150)(1;210)(1;270)(1;330)
\psline(1.932;15)(1;30)(1.932;45) \psline(1.932;75)(1;90)(1.932;105)
\psline(1.932;135)(1;150)(1.932;165) \psline(1.932;195)(1;210)(1.932;225)
\psline(1.932;255)(1;270)(1.932;285) \psline(1.932;315)(1;330)(1.932;345)
\psline(1;30)(1;210) \psline(1;90)(1;270) \psline(1;150)(1;330)
\endpspicture}
\rput(5.25,3.5){\pspicture(-2,-2)(2,2) 
\pspolygon(0,-1.866)(.866,-1.366)(1.366,-.5)(1.366,.5)(.5,1)(-.5,1)
    (-1.366,.5)(-1.366,-.5)(-.866,-1.366)
\pspolygon(0,-.866)(.5,0)(-.5,0)
\psline(1.366,-.5)(.5,0)(.5,1) \psline(-1.366,-.5)(-.5,0)(-.5,1)
\psline(-.866,-1.366)(0,-.866)(.866,-1.366)
\psline(1.366,.5)(.5,0) \psline(-1.366,.5)(-.5,0) \psline(0,-1.866)(0,-.866)
\endpspicture}
\rput(8.5,3){\pspicture(-2,-1.5)(0,2) 
\pspolygon(-0,.5)(-.5,1.366)(-1.366,.866)(-1.866,0)(-1.366,-.866)
    (-.5,-1.366)(-0,-.5)
\psline(-.866,0)(-0,.5) \psline(-.866,0)(0,-.5) \psline(-.866,0)(-1.366,.866)
\psline(-.866,0)(-1.866,0) \psline(-.866,0)(-1.366,-.866)
\endpspicture}
\endpspicture \eatline \end{fullfigure}

If $Q$ is a single cross polytope, then $E$ is a 24-cell. We can also glue
simplices onto subsets of the facets of the cross polytope. The new dihedral
angles formed by such a gluing are $5\pi/6$. The resulting compound is convex
as long as no two glued cross polytope facets share a ridge.  We used a
computer program to list the combinatorially distinct ways of choosing a subset
of nonadjacent facets of the cross polytope; the results may be summarized as
follows. In addition to the 24-cell, this yields 20 new $2$-simple,
$2$-simplicial polytopes.

\begin{proposition}
There are exactly $21$ distinct simplicial edge-tangent compounds composed of
one regular $4$-dimensional cross polytope and $k\ge0$ simplices, according to
the following table:
$$\begin{tabular}{r@{\;}|@{\;}ccccccccc@{\;}|@{\;}c}
$k$&0&1&2&3&4&5&6&7&8&Total\\ \hline
\#&1&1&3&3&6&3&2&1&1&21\\
\end{tabular}$$
\label{p:crosspoly} \end{proposition}

We can also confirm that every square-triangle jewel arises as the edge link of
an edge-tangent $4$-dimensional compound.  For every jewel other than the one
in the center, we can form a convex \emph{edge bouquet} consisting of simplices
and cross polytopes that meet at an edge: we replace each triangle by a simplex
and each square by a cross polytope.  Since the central jewel has two adjacent
squares, its edge bouquet is not convex.  Instead we glue two cross polytopes
along a facet so that the 4 ridges of that facet are flush, \ie, their dihedral
angle is $\pi$.  Thus we can ``caulk'' each such ridge with three simplices
that share the ridge.  The central jewel is the link of 6 of the edges of the
resulting compound of 2 cross polytopes and 12 simplices.

Simplices and regular cross polytopes combine to form many other edge-tangent
simplicial polytopes and hence $2$-simple, $2$-simplicial
polytopes. In particular, these methods lead to the following theorem.

\begin{theorem}
There are infinitely many combinatorially distinct $2$-simple, $2$-simplicial
facet-tangent $4$-polytopes.
\label{th:inf} \end{theorem}

\begin{proof} We glue $n$ cross polytopes end-to-end. Each adjacent pair
produces 4 flush ridges that we caulk with chains of three simplices.  The
facets to which these simplices are glued are not adjacent and so do not
produce any further concavities.
\end{proof}

The chain of $n$ cross polytopes has $f$-vector
\[
(4n+4,18n+6,28n+4,14n+2).
\]
Filling a concavity adds $(2,9,14,7)$, so after filling the $4(n-1)$
concavities we get a simplicial polytope $Q$ with
\[
f(Q)\;=\;(12n- 4,54n-30,84n-54,42-26),
\]
which yields a $2$-simple, $2$-simplicial $4$-polytope $E$
with
\[
f(E)\;=\;(54n-30,252n-156,252n-156,54n-30)
\]
by \prop{p:e-construction}.

\begin{remark} Every $2$-simple, $2$-simplicial $4$-polytopes that we know is
combinatorially equivalent to one which circumscribes the sphere.  Are there
any that are not?

We do know a few $2$-simple, $2$-simplicial $4$-polytopes which are not
$E$-polytopes.  Trivially there is the simplex. There are a few others that
arise by the fact that the $24$-cell is the $E$-polytope of a cross polytope
in 3 different ways.  Color the vertices of a $24$-cell red, green, and blue,
so that the vertices of each color span a cross polytope. If we cap one facet
of a cross-polytope by a simplex and apply the $E$-construction, the result is
a $24$-cell in which 6 facets that meet at 1 vertex are replaced by 10 facets
and 4 vertices.  If the replaced vertex is red, the replacement can be induced
by capping either blue cross polytope or the green cross polytope and then
applying the $E$-construction; the position of the replacement differs between
the two cases. If we replace two different red vertices, one by capping the
green cross polytope and the other by capping the blue cross polytope, then
the resulting polytope is $2$-simple and $2$-simplicial but not an
$E$-polytope.  This construction has several variations: for example, we can
also replace three vertices, one of each color.
\end{remark}

\subsection{Compounds involving the $600$-cell}
\label{s:600}

If $Q$ is the $600$-cell, then $E$ is a $2$-simple, $2$-simplicial polytope
with $f$-vector $(720,3600,3600,720)$ and fatness exactly $5$.  Again, we can
glue simplices onto any subset of nonadjacent facets of the $600$-cell,
creating convex compounds with dihedral angle $14\pi/15$. We did not count the
(large) number of distinct ways of choosing such a subset, analogous to
\prop{p:crosspoly}. It is not possible to glue a cross polytope to a
$600$-cell, because that would create an $11\pi/10$ angle (\ie, a concave
dihedral of $9\pi/10$) which cannot be filled by additional simplices or cross
polytopes.

The large dihedral angles of the $600$-cell make it difficult to form
compounds from it, but we can modify it as follows to create smaller dihedrals.
Remove a vertex and form the convex hull of the remaining 119 vertices.
The resulting convex polytope has 580 of the $600$-cell's tetrahedral facets
and one icosahedral facet. The pentagonal edge link (\fig{f:cutpents}(a)) of
the edges bordering this new facet become modified in a similar way, by
removing one vertex and forming the convex hull of the remaining four vertices
(\fig{f:cutpents}(b)), which results in a trapezoid; thus, the hyperbolic
dihedrals at the ridges around the new facet are $2\pi/5$.

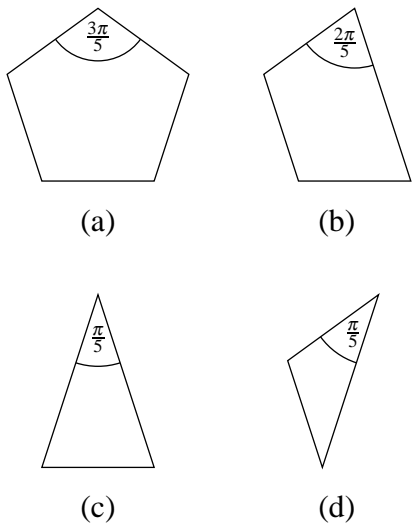
\begin{fullfigure}{f:cutpents}{Edge figures of (a) a $600$-cell, (b) a
    $600$-cell with one vertex removed, (c) a $600$-cell with two removed
    vertices, and (d) an icosahedral cap}
\pspicture(0,-1)(9,10)
\rput(2,8){\pspicture(-2,-2)(2,2)
\pspolygon(2;90)(2;162)(2;234)(2;306)(2;18)
\uput[270](2;90){$\frac{3\pi}{5}$} \psarc(2;90){1.1}{216}{324}
\endpspicture} \rput(2,5.5){\large (a)}
\rput(7,8){\pspicture(-2,-2)(1.25,2)
\pspolygon(2;90)(2;162)(2;234)(2;306)
\uput{.4}[253](2;90){$\frac{2\pi}{5}$} \psarc(2;90){1.25}{216}{288}
\endpspicture} \rput(7,5.5){\large (b)}
\rput(2,2){\pspicture(-1.25,-2)(1.25,2)
\pspolygon(2;90)(2;234)(2;306)
\uput{.65}[270](2;90){$\frac{\pi}{5}$} \psarc(2;90){1.5}{252}{288}
\endpspicture} \rput(2,-.5){\large (c)}
\rput(7,2){\pspicture(-2,-2)(.25,2)
\pspolygon(2;90)(2;162)(2;234)
\uput{.65}[234](2;90){$\frac{\pi}{5}$} \psarc(2;90){1.5}{216}{252}
\endpspicture} \rput(7,-.5){\large (d)}
\endpspicture \eatline \end{fullfigure}

This cut polytope is not simplicial, but we may glue two of these polytopes
together along their icosahedral facets, forming a simplicial polytope with
$4\pi/5$ dihedrals along the glued ridges.  This compound's new edge links are
hexagons formed by gluing pairs of trapezoids (\fig{f:link600}(a)).  The
same cutting and gluing process may be repeated to form a sequence or tree of
$600$-cells, connected along cuts that do not share a ridge.  For such a chain
or tree formed from $n$ cut $600$-cells, the $f$-vector may be computed as
\[
f(Q)\;=\;(106n+14,666n+54,1120n+80,560n+40),
\]
so the $E$-construction yields
\[
f(E)\;=\;(666n+54,3360n+240,3360n+240,666n+54)
\]
and thus a fatness of
\[
\phi(E)\;=\;\frac{3360n+240}{666n+54}\ \
\longrightarrow\ \ \frac{560}{111}\ \approx\ 5.045045.
\]
Thus the fatness of the $2$-simple, $2$-simplicial polytopes formed by such
compounds improves slightly on that formed from the $600$-cell alone.

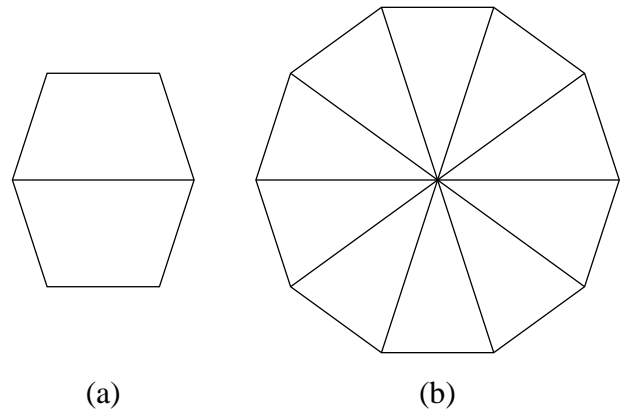
\begin{fullfigure}{f:link600}{Edge links of compounds of cut $600$-cells}
\pspicture(0,-5)(13.5,4)
\rput(2.5,0){\pspicture(-2,-2.5)(2,2.5)
\pspolygon(1.176,2.236)(-1.176,2.236)(-1.902,0)
    (-1.176,-2.236)(1.176,-2.236)(1.902,0)
\psline(-1.902,0)(1.902,0)
\endpspicture}
\rput(2.5,-4.5){\large (a)}
\rput(9.5,0){\pspicture(-4,-4)(4,4)
\pspolygon(3.804;0)(3.804;36)(3.804;72)(3.804;108)(3.804;144)
    (3.804;180)(3.804;216)(3.804;252)(3.804;288)(3.804;324)
\psline(3.804;0)(3.804;180) \psline(3.804;36)(3.804;216)
\psline(3.804;72)(3.804;252) \psline(3.804;108)(3.804;288)
\psline(3.804;144)(3.804;324)
\endpspicture}
\rput(9.5,-4.5){\large (b)}
\endpspicture
\end{fullfigure}

It is also possible to form compounds involving $600$-cells which have been cut
by removing several vertices (as described above) so that two of the resulting
icosahedral facets meet at a ridge. Each edge link at this ridge is an
isosceles triangle formed by removing two vertices from a pentagon
(\fig{f:cutpents}(c)). Thus the dihedral angle of the triangular ridge between
the icosahedral facets is $\pi/5$. We can therefore form compounds in which ten
of these doubly-cut $600$-cells share a triangle, whose edges links are a
regular decagon cut into ten isosceles triangles (\fig{f:link600}(b)).
Yet other compounds of cut $600$-cells and simplices are possible,
although we do not need them here.

The cut $600$-cells also form more complicated compounds which require some
group-theoretic terminology to explain.  The vertices of a regular $600$-cell
form a 120-element group under (rescaled) quaternionic multiplication, the
binary icosahedral group.  This group has a 24-element subgroup, the binary
tetrahedral group, which also arises as the units of the Hurwitz integers (see
Conway and Sloane \cite[\S2.2.6,8.2.1]{CS:splag}).  Let $A$ be the convex hull
of the other 96 vertices of the $600$-cell;  \ie, $A$ is formed by cutting 24
vertices from $600$-cell in the above manner. The resulting polytope is the
``snub$\{3,4,3\}$'' (snub 24-cell) of Coxeter
\cite[\S8.4,8.5]{Coxeter:regpoly}.  Its $f$-vector is $(96,432,480,144)$.
Every icosahedral facet of $A$ is adjacent to $8$ other icosahedron facets, as
well as to $12$ tetrahedra. Thus $A$ has $96$ icosahedron-icosahedron ridges.

We can build new hyperbolic, edge-tangent, simplicial polytopes by gluing
copies of $A$ along icosahedral faces and capping the remaining icosahedral
facets with pyramidal caps of the type that we had cut off to form $A$. (The
edge links of such a cap $C$ are given by \fig{f:cutpents}(d).) The resulting
polytope $Q$ will be convex if at each icosahedral-icosahedral ridge of a copy
of $A$, either $10$ copies of $A$ meet, or two caps and one or two copies of
$A$ do.  Also at each icosahedral-tetrahedral ridge of a copy of $A$, either
two copies of $A$ or one each of $A$ and a cap must meet. If two copies of $A$
meet (in an icosahedral facet $F$), then they differ by a reflection through
$F$. These reflections generate a discrete hyperbolic reflection group $\Gamma$
since the supporting hyperplanes of the icosahedral facets (the facets of a
hyperideal $24$-cell, whose ridges are also ridges of $A$!) satisfy the Coxeter
condition: When they meet, they meet at an angle of $\pi/5$, which divides
$\pi$. Thus the copies of $A$ used in $Q$ are a finite subset $\Sigma$ of the
orbit of $A$ under $\Gamma$.  The set $\Sigma$ determines $Q$.

\begin{remark} That it suffices to consider the dihedral angles of adjacent
facets follows from Poincar\'e's covering-space argument:  Let $P$ be a
spherical, Euclidean, or hyperbolic polytope whose dihedral angles divide
$\pi$.  Let $X$ denote the space in which it lives.  Let $Y$ be the disjoint
union of all copies of $P$ in $X$ in every position, and let $Z$ be the
quotient of $Y$ given by identifying two copies of $P$ along a shared facet.
The space $Z$ is constructed abstractly so that $P$ tiles it.

We claim that $Z$ is a covering space of $X$.  Each $p \in Z$ lies in the
interior of some face $F$ of a copy of $P$.  If $F$ is a copy of $P$ or a
facet, this is elementary; if $F$ is a ridge, it follows from the dihedral
angle condition.  Otherwise it follows by applying the covering-space argument
inductively, replacing $X$ by the link $S$ of $F$ and $P$ by $P \cap S$.

Since $Z$ is a covering space, a connected component of $Z$ is a tiling of $X$
by $P$. See Vinberg \cite{Vinberg:reflection} for a survey of hyperbolic
reflection groups.
\end{remark}

There is no one best choice for $\Sigma$, only a supremal limit.  One
reasonable choice for $\Sigma$ is the corona of a copy of $A$, \ie, $A$
together with the set of all images under $\Gamma$ that meet it, necessarily at
a facet or a ridge. The corona of $A$ is depicted by a simplified (and
therefore erroneous) schematic in \fig{f:corona}; the reader should imagine the
correct, more complicated version.  The schematic uses a chemistry notation in
which each copy of $A$ is represented as an \emph{atom}, each pair of copies
that shares a facet is represented as a \emph{bond}, and ``$7A$'' denotes a
chain of 7 atoms.  To extend the terminology, we call 10 copies of $A$ that
meet at a ridge a \emph{ring}. The schematic is simplified in that the central
copy actually has 24 bonds (not 6), each of the neighbors has 9 bonds (not 3),
and there are 96 rings (not 6).

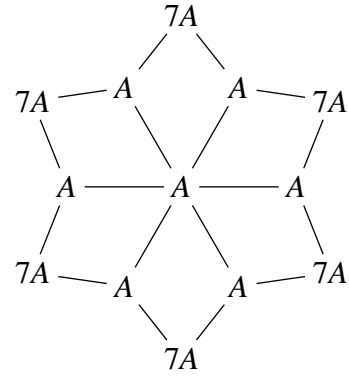
\begin{fullfigure}{f:corona}{An oversimplified schematic of a corona of $A$}
\begin{large} \pspicture(-4,-4)(4,4) \psset{nodesep=3pt}
\rput(0,0){\rnode{a}{$A$}}
\rput(2.4;0){\rnode{b1}{$A$}} \rput(2.4;60){\rnode{b2}{$A$}}
\rput(2.4;120){\rnode{b3}{$A$}} \rput(2.4;180){\rnode{b4}{$A$}}
\rput(2.4;240){\rnode{b5}{$A$}} \rput(2.4;300){\rnode{b6}{$A$}}
\rput(3.6;30){\rnode{c1}{$7A$}} \rput(3.6;90){\rnode{c2}{$7A$}}
\rput(3.6;150){\rnode{c3}{$7A$}} \rput(3.6;210){\rnode{c4}{$7A$}}
\rput(3.6;270){\rnode{c5}{$7A$}} \rput(3.6;330){\rnode{c6}{$7A$}}
\ncline{a}{b1} \ncline{a}{b2} \ncline{a}{b3}
\ncline{a}{b4} \ncline{a}{b5} \ncline{a}{b6}
\ncline{b1}{c1} \ncline{b1}{c6} \ncline{b2}{c2} \ncline{b2}{c1}
\ncline{b3}{c3} \ncline{b3}{c2} \ncline{b4}{c4} \ncline{b4}{c3}
\ncline{b5}{c5} \ncline{b5}{c4} \ncline{b6}{c6} \ncline{b6}{c5}
\endpspicture \end{large} \eatline \end{fullfigure}

Let $Q$ be the union of these copies of $A$ with the remaining icosahedral
facets capped.  To compute the $f$-vector of $Q$ it is easier to view each atom
as a copy of a 600-cell $B$, minus two caps for each bond.  There are $1+24+96
\cdot 7 = 697$ atoms and $24+ 96 \cdot 8 = 792$ bonds in total. Thus $Q$ has
\[
f_3(Q) = 697 f_3(B) - 792 \cdot 2 f_3'(C) = 386520
\]
facets, where $f_3'(C)=30$ is the number of simplicial facets of a cap $C$.

Counting vertices is more complicated.  Let $I$ be an icosahedron
and let $T$ be a triangle.  Then
\begin{align*}
f_0(Q) &= 697 f_0(B) - 792 \cdot 2 f_0(C) + 792 f_0(I) + 96 f_0(T) \\
&= 72840,
\end{align*}
because after the vertices of the caps are subtracted, the vertices of each
icosahedral facet at a bond are undercounted once, and after these are restored
the vertices of each triangle at the center of a ring are undercounted once.
The rest of the $f$-vector of $Q$ follows from the Dehn-Sommerville equations:
\begin{align*}
f_2(Q) &= 2f_3(Q) = 773040, \\
f_1(Q) &= f_0(Q) + f_3(Q) = 459360.
\end{align*}
The polytope $Q$ yields an $E$-polytope with fatness
\[
\phi(E) = 6\,\frac{f_3(Q)}{f_1(Q)} = \frac{3221}{638} \approx 5.048589.
\]
Note that since this bound arises from a specific choice of $\Sigma$
rather than a supremal limit, this is not optimal as a lower
bound of supremal fatness.

\subsection{Kissing numbers}
\label{s:kissing}

As mentioned above, another use of ridge-tangent polytopes $P$ is the average
kissing number problem \cite{Kuperberg:kissing}.  Let $X$ be a finite packing
of (not necessarily congruent) spheres $S^3$, which is equivalent to a finite
sphere packing in $\R^3$ by stereographic projection. The question is to
maximize the average number of kissing points of the spheres in $X$.  If $P$ is
ridge-tangent, its facets intersect the unit sphere $S^3$ in a sphere packing
$X$, in which the spheres kiss at the tangency points of $P$. Thus the
ridge-facet ratio of $P$ is exactly half the average kissing number of $X$.
(Not all sphere packings come from ridge-tangent polytopes in this way.)

The sphere packings due to Kuperberg and Schramm \cite{Kuperberg:kissing} can
be viewed as coming from a compound consisting of a chain or tree of $n$ cut
$600$-cells (\ie, atoms in the sense of \fig{f:corona}).  Their average kissing
numbers are
\[
\kappa\; =\; 2\,\frac{666n+54}{106n+14}\ \
\longrightarrow\ \ \frac{666}{53}\ \approx\ 12.56603.
\]
By contrast if $Q$ is the compound formed from a corona of $A$ as in
\sec{s:600}, then the average kissing number of the corresponding sphere
packing is
\[
\kappa(Q) = 2\,\frac{f_1(Q)}{f_0(Q)} = \frac{7656}{607} \approx 12.61285.
\]
Like the bound on fatness, it is not optimal as a lower bound
on the supremal average kissing number.

Here we offer no improvement on the upper bound
\[
\kappa < 8+4\sqrt{3} \approx 14.92820
\]
from \cite{Kuperberg:kissing}, even though it cannot be optimal either.

\section{3-Spheres}
\label{s:spheres}

In this section we construct a family of strongly regular cellulations of the
$3$-sphere with unbounded fatness.  Indeed, we provide an efficient version of
the construction, in the sense that it requires only polynomially many cells to
achieve a given fatness.  (The construction is also a polynomially effective
randomized algorithm.)  Given $N$, we find a strongly regular cellulation of
$S^3$ with $O(N^{12})$ cells and fatness at least $N$.  Note that there are
also power-law upper bounds on fatness:  An $O(N^{1/3})$ upper bound for the
fatness of a convex $4$-polytope with $N$ vertices follows from work by
Edelsbrunner and Sharir \cite{ES:hyperplane}.  The K\H{o}vari-S\'os-Tur\'an
theorem \cite{KST:problem} (see also \cite[p.~1239]{Bollobas:extremal} and
\cite[Thm.~9.6,p.~121]{PA:combin}) implies an $O(N^{2/3})$ upper bound on the
fatness of strongly regular cellulations of $S^3$, since the vertex-facet
(atom-coatom) incidence graph has no $K_{3,3}$-subgraph, and thus has at most
$O\big(f_0(f_0+f_3)^{2/3}\big)$ edges. Our construction provides an
$\Omega(N^{1/12})$ lower bound.

The inefficient construction is a simpler version which we describe first.  For
every $g>0$, $S^3$ can be realized as
\[
H_1 \cup (S_g \times I) \cup H_2,
\]
a thickened surface of genus $g$ capped on both ends with handlebodies. (This
is obtained from a neighborhood of the standard [unknotted] smooth embedding
of $S_g$ into $S^3$.) If for some $g$ we can find a fat cellulation of $S_g$,
we can realize $S^3$ as a ``fat sausage with lean ends,'' as shown in
\fig{f:sausage}.  We cross the fat cellulation of $S_g$ with an interval
divided into $N$ segments to produce a fat cellulation of $S_g \times I$.  Then
we fix arbitrary strongly regular cellulations of the handlebodies $H_1$ and
$H_2$.  If we make the sausage long enough, \ie, if we take $N \to \infty$, the
fatness of the sausage converges to the fatness of its middle regardless of the
structure of its ends.

\begin{fullfigure}{f:sausage}{$S^3$ as a fat sausage with lean ends}
\begin{large} \psset{unit=.2in} \pspicture(-8,-2.5)(8,3)
\multips(-6,0)(4,0){4}{
    \psbezier(.7,-1)(.7,-1.5)(.35,-2)(0,-2)
    \psbezier(.7,-1)(.7,-.5)(.35,-.5)(.35,0)
    \psbezier(.35,0)(.35,.5)(.7,.5)(.7,1)
    \psbezier(.7,1)(.7,1.5)(.35,2)(0,2)
    \psbezier(-.7,-1)(-.7,-1.5)(-.35,-2)(0,-2)
    \psbezier(-.7,-1)(-.7,-.5)(-.35,-.5)(-.35,0)
    \psbezier(-.35,0)(-.35,.5)(-.7,.5)(-.7,1)
    \psbezier(-.7,1)(-.7,1.5)(-.35,2)(0,2)
    \psbezier(0,.65)(.25,.9)(.25,1.15)(0,1.4)
    \psbezier(0,.65)(-.25,.9)(-.25,1.15)(0,1.4)
    \psbezier(0,-.65)(.25,-.9)(.25,-1.15)(0,-1.4)
    \psbezier(0,-.65)(-.25,-.9)(-.25,-1.15)(0,-1.4)}
\psline(-6,2)(6,2)\psline(-6,-2)(6,-2)
\psarc(6,0){2}{270}{90}\psarc(-6,0){2}{90}{270}
\psline[linestyle=dashed](-6,1.4)(6,1.4)
\psline[linestyle=dashed](-6,.65)(6,.65)
\psline[linestyle=dashed](-6,-1.4)(6,-1.4)
\psline[linestyle=dashed](-6,-.65)(6,-.65)
\rput(-7.25,0){$H_1$}\rput(0,0){$S_g \times I$}\rput(7.25,0){$H_2$}
\endpspicture \end{large} \eatline \end{fullfigure}

It remains only to show that there are strongly regular fat cellulations of
surfaces.  The surface $S_g$ has perfect cellulations with $f$-vector
$(1,2g,1)$.  Such a cellulation is obtained by gluing pairs of sides of a
$4g$-gon in such a way that all vertices are identified. It has fatness $g$,
and it exists for arbitrarily large $g$, but it is far from regular. However,
its lift to the universal cover $\tilde{S}_g$ is strongly regular, since any
such cellulation can be represented by a tiling of the hyperbolic plane by
\emph{convex} polygons.  (Indeed, if we take the regular $4g$-gon with angles
of $\pi/(2g)$, which is certainly convex, then its edges and angles are
compatible with any perfect cellulation.)  Moreover, Mal'cev's theorem
\cite{Malcev:faithful}, states that finitely generated matrix groups are
residually finite; this implies that every closed hyperbolic manifold admits
intermediate finite covers with arbitrarily large injectivity radius.  (See
\cite[\S4]{Kuperberg:saturated} for a detailed exposition.) In particular $S_g$
admits an intermediate cover $\hS_g$ whose injectivity radius exceeds the
diameter of a $2$-cell.  The cellulation of $\hS_g$ is then strongly
regular.  Its genus is much larger than $g$, but its fatness is still $g$.

The efficient construction is the same: It only requires careful choices for
the finite cover $\hS_g$ and for the handlebodies $H_1$ and $H_2$. Among
the perfect cellulations of $S_g$, a convenient one for us is a $4g$-gon with
opposite edges identified.  We describe the fundamental group $\pi_1(S_g)$
using this cellulation.  As shown in \fig{f:sg}, we number the edges
$x_0,\dots,x_{4g-1}$ consecutively, so that
\eq{e:ident}{x_i^{} = x_{2g+i}^{-1}}
and
\eq{e:cell}{x_0x_1\cdots x_{4g-1} = 1.}
(We interpret the indices as elements of $\Z/(4g)$.)
Equation~\eqref{e:ident} expresses the identifications, while
equation~\eqref{e:cell} expresses the boundary of the $2$-cell.

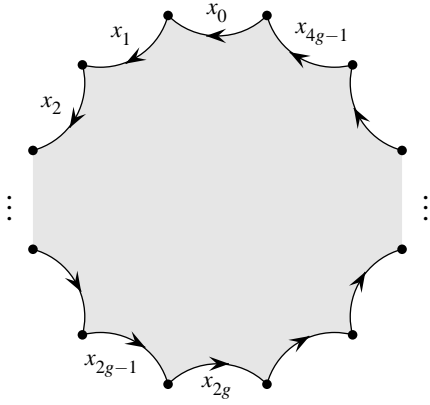
\begin{fullfigure}{f:sg}{Labelling edges of $S_g$}
\pspicture(-4.5,-4.5)(4,4.5) \sgpie
\uput[90](3.435;90){$x_0$} \uput[120](3.435;120){$x_1$}
\uput[150](3.435;150){$x_2$} \uput[240](3.435;240){$x_{2g-1}$}
\uput[270](3.435;270){$x_{2g}$} \uput[60](3.435;60){$x_{4g-1}$}
\eatline \endpspicture \end{fullfigure}

We construct $\hS_g$ as a tower of two abelian finite covers, which
together form an irregular cover of $S_g$. (No abelian cover of $S_g$ is
strongly regular.) The surface $S_g$ satisfies the usual requirements of
covering-space
theory (see Fulton \cite[\S\ 13b, 14a]{Fulton:gtm}): It is a connected,
locally path-connected, and locally simply connected space. For any finite
group $A$, the \emph{$A$-coverings} of $S_g$ are covering spaces of the form
$Y$ with $Y/A=S_g$, where $A$ acts properly discontinuously on $Y$. These
coverings, up to isomorphism, are classified by the set of group homomorphisms
$\Hom(\pi_1(S_g,x),A)$ \cite[Thm.~14.a]{Fulton:gtm}. Furthermore, if $A$ is
abelian, then every such homomorphism maps the commutators in $\pi_1(S_g,x)$ to
zero, so the $A$-coverings are classified by $\Hom(H_1(S_g,\Z),A)$. In other
words, if $A$ is any abelian group with $n$ elements, then every homomorphism
$\sigma: H_1(S_g)\rightarrow A$ defines an $n$-fold abelian covering of $S_g$.
(If the homomorphism is not surjective, then the covering space is not
connected \cite[p.~193]{Fulton:gtm}. In this case we use a connected component
of the covering space, which has the same fatness but smaller genus.)

Now assume that $q = 4g+1$ is a prime power and let $\alpha$ generate the
cyclic group $\F_q^*$, where $\F_q$ is the field with $q$ elements.  We can
fulfill the assumption by changing $g$ by a bounded factor.  (Most simply we
can let $q = 5^k$.  Or we can let $q$ be prime, so that $\F_q=\Z/q$, by a form
of Bertrand's postulate for primes in congruence classes.  This result dates to
the 19th century; see Erd\H{o}s \cite{Erdos:primzahlen} for an elementary
proof.) Define a homomorphism $\sigma:H_1(S_g) \to \F_q$ by $\sigma([x_i]) =
\alpha^i$, where $[x_i]$ is the 1-cycle (or homology class) represented by the
loop $x_i$.  Since $\alpha^{2g} = -1$, the definition of $\sigma$ is consistent
with equation~\eqref{e:ident}. Consistency with \eqref{e:cell} is then
automatic. Let $S_g'$ be the finite cover corresponding to $\sigma$.

To prepare for the subsequent analysis of $\hS_g$, we give an explicit
combinatorial description of $S'_g$.  Let $F^0$ be a $2$-cell of $S'_g$
and label its vertices
$$v_0^0,v_1^0,\ldots,v_{4g-1}^0$$
in cyclic order.  See \fig{f:sg1}.

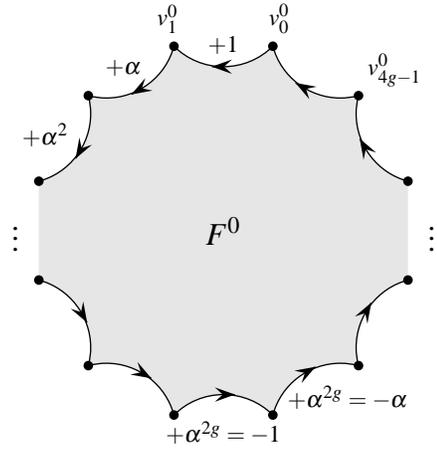
\begin{fullfigure}{f:sg1}{Labelling the vertices of $F^0$}
\pspicture(-5,-4.75)(5,4.75) \sgpie \rput(0,0){\large $F^0$}
\uput[90](3.435;90){$+1$} \uput[120](3.435;120){$+\alpha$}
\uput[150](3.435;150){$+\alpha^2$}
\uput{10pt}[270](3.435;270){$+\alpha^{2g} = -1$}
\uput[300](3.435;300){$+\alpha^{2g} = -\alpha$}
\uput[45](4;45){$v^0_{4g-1}$} \uput[75](4;75){$v^0_0$}
\uput[105](4;105){$v^0_1$}
\eatline \endpspicture \end{fullfigure}

For each $s \in \F_q$, let $F^s$ and $v_k^s$ be the images of $F^0$ and $v_k^0$
under the action of $s$. Since $S_g$ has only one vertex, the vertices of
$S_g'$ may be identified with $\F_q$. Thus, if we identify $v_0^0$ with
$0\in\F_q$, then the action of $\F_q$ will identify $v_0^s$ with $s\in\F_q$.
The structure of $\sigma$ further implies that
\[
v_k^s\ =\ s+1+\alpha+\alpha^2+\ldots+\alpha^{k-1}\ =\
    s+\frac{\alpha^k-1}{\alpha-1}
\]
for all $k\in\Z/(4g)$ and $s\in\F_q$.  See \fig{f:sg2}.

\begin{fullfigure}{f:sg2}{Labelling the vertices of $F^s$}
\pspicture(-7,-4.75)(5.5,4.75) \sgpie \rput(0,0){\large $F^s$}
\uput[75](4;75){$v^s_0 = s$}
\uput[105](4;105){$v^s_1 = s+1$}
\uput[120](4;135){$v^s_2 = s+1+\alpha$}
\uput[135](4;165){$v^s_3 = s+\tfrac{\alpha^3-1}{\alpha-1}$}
\eatline \endpspicture \end{fullfigure}

Using this explicit description, it is routine to verify the following
(remarkable) properties of the surface $S'_g$.

\begin{lemma}
The cellulation of the abelian cover $S'_g$ is regular and has $f$-vector
\[
f(S'_g) = (q,2gq,q)=(q,\binom{q}{2},q).
\]
Every facet has $q-1=4g$ vertices, while every vertex has degree $q-1=4g$. The
graph (or $1$-skeleton) of $S_g'$ is the complete graph on $q+1$ vertices. The
dual graph is also complete; any two facets share exactly one edge (as well as
$q-4$ other vertices).
\label{l:regular} \end{lemma}

\begin{proof}
In view of the combinatorial description above (\fig{f:sg2}), all these
facts follow from simple computations in the field $\F_q$:
\begin{itemize}\setlength{\itemsep}{0pt}
\item $S_g'$ is regular --- for each $s\in\F_q$, the vertex labels
$s+\frac{\alpha^k-1}{\alpha-1}$ ($0\le k<4g$) are distinct.
\item The $1$-skeleton of $S_g'$ is complete --- for $v,v'\in\F_q$, $v\neq v'$
there is a unique $s\in\F_q$ and $k\in\Z/(4g)$ with
\[
v=s+\frac{\alpha^k-1}{\alpha-1}\qquad
    v'=s+\frac{\alpha^{k+1}-1}{\alpha-1}.
\]
\item The dual graph of $S_g'$ is complete --- for $s,s'\in\F_q$, $s\neq s'$,
there are unique $k,\ell\in\Z/(4g)$ such that
\begin{align*}
s+\frac{\alpha^k-1}{\alpha-1}&=s'+\frac{\alpha^{\ell+1}-1}{\alpha-1} \\
s+\frac{\alpha^{k+1}-1}{\alpha-1}&=s'+\frac{\alpha^\ell-1}{\alpha-1}.
\end{align*}
\end{itemize}
\end{proof}

\begin{theorem}
Let $n \ge 128 g^4$, and let
\[
\rho:H_1(S_g')\ \to\ \Z/n
\]
be a randomly chosen homomorphism, and let $\hS_g$ be the finite cover of
$S_g'$ corresponding to $\rho$. Then with probability more than $\frac12$, the
cellulation of\/ $\hS_g$ is strongly regular.
\label{th:prob}
\end{theorem}

In order to prove \thm{th:prob}, we need to more explicitly describe the
condition of strong regularity as it applies to $\hS_g$. Let $X$ be a regular
cell complex and suppose that its universal cover $\tilde{X}$ is strongly
regular.  Recall that the \emph{star} $\st(v)$ of a vertex $v$ in $X$ is the
subcomplex generated by the cells that contain $v$. The complex $X$ is strongly
regular if and only if the star of each vertex is.  Suppose that $\tilde{v} \in
\tilde{X}$ projects to $v\in X$. Then the star $\st(\tilde{v})$, which is
strongly regular, projects to the star $\st(v)$.  The latter is strongly
regular if and only if the projection is injective.  In other words, $X$ is
strongly regular if and only if the stars of $\tilde{X}$ embed in $X$. If $X$
is not strongly regular, then there must be a path $\tilde{\ell}$ connecting
distinct vertices of $\st(\tilde{v})$ which projects to a loop $\ell$ in
$\st(v)$.  We say that such a loop \emph{obstructs strong regularity}.  We can
assume that $\tilde{\ell}$ is a pair of segments properly embedded in distinct
cells in $\st(\tilde v)$, with only the end-points of the segments on the
boundary of the cells, which implies that $\ell$ is embedded if $X$ is regular.
\fig{f:obstruct} gives an example of such a loop $\ell$ in a regular
cellulation of a torus.

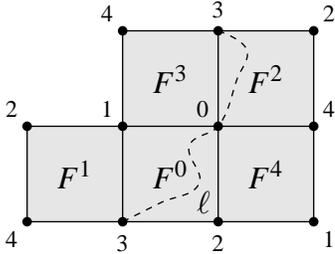
\begin{fullfigure}{f:obstruct}{A loop $\ell$ that obstructs strong regularity
    in the torus $S_1'$}
\pspicture(0,0)(6,4)
\pspolygon[linestyle=none,fillstyle=solid,fillcolor=gray90]
    (0,0)(0,2)(2,2)(2,4)(6,4)(6,0)
\psline(0,0)(6,0) \psline(0,2)(6,2) \psline(2,4)(6,4) \psline(0,0)(0,2)
\psline(2,0)(2,4) \psline(4,0)(4,4) \psline(6,0)(6,4)
\multips(0,0)(2,0){4}{\qdisk(0,0){.1}} \multips(0,2)(2,0){4}{\qdisk(0,0){.1}}
\multips(2,4)(2,0){3}{\qdisk(0,0){.1}}
\rput(1,1){\large $F^1$} \rput(3,1){\large $F^0$} \rput(5,1){\large $F^4$}
\rput(3,3){\large $F^3$} \rput(5,3){\large $F^2$}
\uput[225](0,0){$4$} \uput[270](2,0){$3$} \uput[270](4,0){$2$}
\uput[315](6,0){$1$}
\uput[135](0,2){$2$} \uput[135](2,2){$1$} \uput[135](4,2){$0$}
\uput[45](6,2){$4$}
\uput[135](2,4){$4$} \uput[90](4,4){$3$} \uput[45](6,4){$2$}
\pscurve[linestyle=dashed](2,0)(2.8,.4)(3.6,.8)(3.4,1.6)(4,2)(4.4,3)
    (4.6,3.6)(4,4)
\rput(3.7,.4){\large $\ell$}
\endpspicture
\end{fullfigure}

In our case, the surfaces $S_g'$ are regular, but they have many obstructing
loops.  \thm{th:prob} asserts that, with non-zero probability, all such loops
lengthen when lifted to $\hS_g$.

\begin{lemma} No loops in $S_g'$ that obstruct strong regularity are
null-homologous.  Furthermore, all obstructing loops represent indivisible
elements in $H_1(S_g')$.
\label{l:notnull} \end{lemma}

\begin{proof} In brief, they are indivisible because they are embedded,
and they are too short to be null-homologous.

By \lem{l:regular}, $S_g'$ is regular.  By the discussion after the statement
of \thm{th:prob}, each obstructing loop $\ell$ is embedded. If $\ell$ separates
$S_g'$, then it is null-homologous. If $\ell$ does not separate $S_g'$, then it
is indivisible in homology.  (To show this, we can appeal to the classification
of surfaces by cutting $S_g'$ along $\ell$.  The classification implies that
all non-separating positions for $\ell$ are equivalent up to homeomorphism of
$S_g'$.  It is easy to find a standard position for $\ell$ in which it is
indivisible in homology.) Thus it remains to show that no obstructing loop is
null-homologous.

First, we claim that any obstructing loop $\ell$ can be supported on fewer than
$4g$ edges of the $1$-skeleton of $S_g'$.  We homotop the two segments of
$\ell$ to the boundaries of the 2-cell containing them, giving them each at
most $2g$ edges.  Thus $\ell$ is represented by a sequence of at most $4g$
edges in $S_g'$. The case of exactly $4g$ edges does not occur, since the
endpoints of the loop coincide, and no two vertices $v\neq v'$ of $S_g'$ are
opposite vertices in two different facets $F^s$.  As in \lem{l:regular}, this
follows from the fact that for $v,v'\in\F_g$, $v\neq v'$, there are unique
$s\in\F_q$ and $k\in\Z/(4g)$ with
\[
v = s+\frac{\alpha^k-1}{\alpha-1}\qquad
v' = s+\frac{\alpha^{k+2g}-1}{\alpha-1}.
\]

Second,
we claim that any null-homologous loop in the $1$-skeleton
of $S_g'$ contains at least $4g$ edges.  In other we if $f$ is a $2$-chain on
$S_g'$ and $\d f \ne 0$, then $|\d f| \ge 4g$. Since $S_g'$ is orientable, we
can regard $f$ as a function on its $2$-cells.  Since $f$ is non-constant, it
attains some value $t$ on $k$ $2$-cells with $0 < k < q$.  Since any two
$2$-cells share an edge by \lem{l:regular}, these $2$-cells share
\[
k(q-k) \ge 4g
\]
edges with the complementary set of $2$-cells, of which there are $q-k$.  Since
$\d f$ is non-zero on these edges, $|\d f| \ge 4g$, as desired.
\end{proof}

\begin{proof}[Proof of \thm{th:prob}] In brief, $S_g'$ has fewer than $64g^4$
obstructing loops $\ell$.  For each one,
\[
P\bigl[\rho([\ell]) = 0\bigr] = \frac1n.
\]
The expected number of obstructing loops that lift from $S_g'$ to $\hS_g$
without lengthening is less than $64g^4/n \le \frac12$.  Thus there is a good
chance that all obstructing loops lengthen.

The homology group $H_1(S'_g)$ is a finitely generated free abelian group: It
is isomorphic to $\Z^d$, for $d=1+q(g-1)$. Thus it admits $n^d$ homomorphisms
$\rho$ to $\Z/n$.  Since this is a finite number, choosing one uniformly at
random is well-defined.  If $c$ is any indivisible vector in $\Z^d$, then it is
contained in a basis, and thus $\rho(c)$ is equidistributed. In particular, if
$\ell$ is an obstructing loop, then $[\ell]$ is indivisible by \lem{l:notnull},
so $\rho([\ell])$ is equidistributed in $\Z/n$.

It remains only to bound the number of obstructing loops in $S_g'$. A star
$\st(v)$ in $\tilde{S}_g$ has $4g(4g-2)$ points other than $v$ itself.
Without loss of generality, $v$ projects to $0$ in $S_g'$.  In this case the
other vertices are equidistributed among the $4g$ non-zero elements of
$\F_q$.  Therefore $\st(v)$ has
\[
4g\binom{4g-2}{2}\ =\ \frac{4g(4g-2)(4g-3)}{2}
\]
pairs of arcs connecting $v$ to two vertices that are the same in $S_g'$.
These pairs represent all two-segment obstructing loops that pass through $0$
and a nonzero vertex $v'$ (and some of these loops are homotopic).  If we count
such pairs of arcs for any pair of distinct vertices $v,v'$
of $S_g'$, then we find that the total number is not more than
\[
\binom{4g+1}{2}\binom{4g-2}{2}\ =\
\frac{(4g+1)4g(4g-2)(4g-3)}{4} < 64g^4,
\]
as desired.
\end{proof}

\begin{question} For each $g>1$, what is the maximum fatness of a strongly
regular cellulation of a surface of genus $g$? Equivalently, how many edges are
needed for a strongly regular cellulation of a surface of genus $g$?
\end{question}

\begin{remark}
One interesting alternative to the construction of $S_g'$ is to assume instead
that $q = 4g-1$ is a prime power, and to let $\alpha$ be an element of order
$4g$ in $\F_{q^2}$.  The resulting $q^2$-fold cover $S_g''$ is almost strongly
regular: the only obstructing loops are those that are null homologous in
$S_g$.  Another interesting surface is the modular curve $X(2p)$, where $p$ is
a prime \cite[\S13]{Silverman:gtm}. The inclusion $\Gamma(2p) \subset
\Gamma(2)$ of modular groups induces a projection from $X(2p)$ to the modular
curve $X(2)$, which is a sphere with three cusp points. If we connect two of
these points by an arc which avoids the third, it lifts to a cellulation of
$X(2p)$ with $f$-vector $(p^2-1,\frac{p(p^2-1)}{2},\frac{p^2-1}{2})$. Like
$S_g''$, it has a few obstructing loops. Unfortunately we do not know a way to
use either $S_g''$ or $X(2p)$ to make fat surfaces of lower genus (or
equivalently fewer cells) than $\hS_g$.
\end{remark}

Since \thm{th:prob} provides us with efficient fat surfaces $\hS_g$, the
construction of fat cellulations of $S^3$ only requires efficient cellulations
of the handlebodies $H_1$ and $H_2$ and an efficient way to attach them to
$\hS_g$. In our construction the handlebody cellulations are \emph{a priori}
unrelated to the cellulation of $\hS_g$.  Rather they are transverse after
attachment, and each point of intersection will become a new vertex.  Thus the
question is to position the cellulations to minimize their intersection.

We describe the cellulations in three stages:  first, a
dissection of $H_1$ and $H_2$ individually into $3$-cells; second, their
relative position; and third, their position relative to the cellulation of
$\hS_g$.

Let $\hg$ be the genus of $\hS_g$.  A handlebody $H$ of genus $\hg$ can
be formed by identifying $\hg$ pairs of disks on the surface of a $3$-cell.
The result is a dissection of $H$ into $\hg$ $2$-cells and one $3$-cell,
although it is not a cell complex because there are no $1$-cells or $0$-cells.
We can still ask whether such a dissection is regular or strongly regular; this
one is neither.  However, if we replace each $2$-cell by 3 parallel $2$-cells,
it becomes strongly regular.  An example of the resulting dissection $A$ is
shown in \fig{f:tripled}.

\begin{fullfigure}{f:tripled}
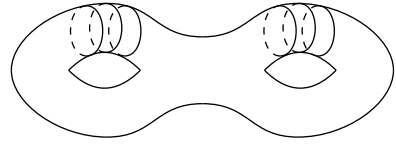
    {A strongly regular dissection $A$ of a handlebody of genus 2}
\pspicture(-4,-2)(4,3)
\psbezier(-2,1.4)(-3,1.4)(-4,0.7)(-4,0)
\psbezier(-2,1.4)(-1,1.4)(-1,0.7)(0,0.7)
\psbezier(0,0.7)(1,0.7)(1,1.4)(2,1.4)
\psbezier(2,1.4)(3,1.4)(4,0.7)(4,0)
\psbezier(-2,-1.4)(-3,-1.4)(-4,-0.7)(-4,0)
\psbezier(-2,-1.4)(-1,-1.4)(-1,-0.7)(0,-0.7)
\psbezier(0,-0.7)(1,-0.7)(1,-1.4)(2,-1.4)
\psbezier(2,-1.4)(3,-1.4)(4,-0.7)(4,0)
\psbezier(1.3,0)(1.8,0.5)(2.3,0.5)(2.8,0)
\psbezier(1.3,0)(1.8,-0.5)(2.3,-0.5)(2.8,0)
\psbezier(-1.3,0)(-1.8,0.5)(-2.3,0.5)(-2.8,0)
\psbezier(-1.3,0)(-1.8,-0.5)(-2.3,-0.5)(-2.8,0)
\psbezier(-2.025,.375)(-1.575,.375)(-1.575,1.4)(-2.025,1.4)
\psbezier[linestyle=dashed](-2.025,.375)(-2.475,.375)(-2.475,1.4)(-2.025,1.4)
\psbezier(-2.425,.3)(-1.975,.3)(-1.975,1.35)(-2.425,1.35)
\psbezier[linestyle=dashed](-2.425,.3)(-2.875,.3)(-2.875,1.35)(-2.425,1.35)
\psbezier(-1.625,.3)(-1.175,.3)(-1.175,1.35)(-1.625,1.35)
\psbezier[linestyle=dashed](-1.625,.3)(-2.075,.3)(-2.075,1.35)(-1.625,1.35)
\psbezier[linestyle=dashed](2.025,.375)(1.575,.375)(1.575,1.4)(2.025,1.4)
\psbezier(2.025,.375)(2.475,.375)(2.475,1.4)(2.025,1.4)
\psbezier[linestyle=dashed](2.425,.3)(1.975,.3)(1.975,1.35)(2.425,1.35)
\psbezier(2.425,.3)(2.875,.3)(2.875,1.35)(2.425,1.35)
\psbezier[linestyle=dashed](1.625,.3)(1.175,.3)(1.175,1.35)(1.625,1.35)
\psbezier(1.625,.3)(2.075,.3)(2.075,1.35)(1.625,1.35)
\endpspicture \end{fullfigure}

The surface $S_\hg$ (which for our choice of $\hg$ is isomorphic to
$\hS_g$) has another standard perfect cellulation called a \emph{canonical
schema} in the computer science literature \cite{VY:surfaces}.  Using the
labelling in \fig{f:sg}, we identify $x_{4k}^{}$ with $x_{4k+2}^{-1}$ (the even
loops), and $x_{4k+1}^{}$ with $x_{4k+3}^{-1}$ (the odd loops), for all $0 \le
k \le \hg$. By rounding corners we can make each even loop intersect one odd
loop once and eliminate all other intersections between loops. If we interpret
this pattern of loops as the standard Heegaard diagram for $S^3$
\cite{Stillwell:gtm}, then the even loops bound disks in $H_1$ and the odd
loops bound disks in $H_2$.  We can then put in two copies $A_1$ and $A_2$ of
the cell division $A$ so that its $2$-cells run parallel to these loops.

We would like to position $A_1$ and $A_2$ to minimize their intersection with
the cellulation of $\hS_g$. (Note that $A_1$ and $A_2$ do not intersect each
other since $\hS_g \times I$ lies in between.) To this end Vegter and Yap
\cite{VY:surfaces} proved that any cellulation of $S_\hg$ with $n$ edges admits
a position of the canonical schema cellulation with $O(n\hg)$ intersections.
(Strictly speaking the theorem applies to triangulations, but any regular
cellulation of a surface with $n$ edges can be refined to a triangulation with
less than $3n$ edges.) In our case
\[
n = \Theta(\hg) = \Theta(g^6).
\]
Thus, by tripling the edges of the canonical schema in the Vegter-Yap
construction, we can position $A_1$ and $A_2$ so that the lean ends in the
sausage have $O(g^{12})$ vertices.  If we give the fat part of the sausage $N =
g^{12}$ slices, the total $f$-vector of the cellulation of $S^3$ is then
\[
(\Theta(g^{12}),\Theta(g^{13}),\Theta(g^{13}),\Theta(g^{12})),
\]
and its fatness is $\Theta(g)$.  This completes the efficient construction
with unbounded fatness.


\providecommand{\bysame}{\leavevmode\hbox to3em{\hrulefill}\thinspace}

\end{document}